\documentclass[12pt]{article}
\usepackage{amsfonts}
\usepackage{amsmath}
\usepackage{amssymb}

\newtheorem{thm}{Theorem}
\newtheorem{cor}[thm]{Corollary}
\newtheorem{pro}[thm]{Propostion}

\newtheorem{defn}[thm]{Definition}
\newtheorem{rmk}[thm]{Remark}

\newtheorem{con}[thm]{Conjecture}
\newtheorem{exmp}[thm]{Example} 

\def\N{\mathbb{N}}           

\def\P{\mathbb{P}}

\def\C{\mathbb{C}}

\def\ep{\varepsilon}
\newcommand{\beq}{\begin{equation}}\newcommand{\eeq}{\end{equation}}

\newtheorem{prf}{Proof}

\def\N{\mathbb{N}}           

\def\P{\mathbb{P}}

\def\C{\mathbb{C}}

\def\nl{\newline}
\def\ep{\varepsilon}
\def\bo{\nl\phantom{a}\hfill $\Box$\nl}

\usepackage{enumitem}

\usepackage[affil-it]{authblk} 
\usepackage{etoolbox}
\usepackage{lmodern}

\makeatletter
\patchcmd{\@maketitle}{\LARGE \@title}{\fontsize{14}{19.2}\selectfont\@title}{}{}
\makeatother

\usepackage{authblk}

\date{}

\title{\bf Multiplicative functions with sum zero}
\author[1,2]{Ammar Ali Neamah}

\affil[1]{Department of Mathematics, University of Reading, Whiteknights, Reading, UK. }
\affil[2]{ Faculty of Computer Science and Mathematics, University of Kufa,  Najaf, Iraq.}

\begin{document}
  \maketitle

\begin{abstract}  
$CMO$ functions are completely multiplicative functions $f$ for which $\sum_{n=1}^\infty f(n)$ $=0$. These functions were first introduced and studied by Kahane and Sa\"{i}as \cite{KS1}. The main purpose of this paper is to generalise such functions to \mbox{multiplicative} \mbox{functions} and we shall call them $MO$ functions. More \mbox{precisely}, we define $MO$ functions to be multiplicative functions for which $\sum_{n=1}^\infty f(n) = 0$ and $\sum_{k=0}^\infty f(p^k)$  $\ne 0$ for all $p \in \P$. 
We give some properties and find examples of $MO$ functions, as well as pointing out the connection between these functions and the Riemann hypothesis at the end of the paper.

\ 

\noindent
{\em 2010 AMS Mathematics Subject Classification}: 11N64, 11M41, 11N99  \nl
{\em Keywords and phrases}: Multiplicative functions, Riemann Hypothesis
\vspace{0.3in}\end{abstract}

\section{Introduction}
An arithmetical function $f: \N \longrightarrow \C$ is called {\em multiplicative} if $f(1) = 1$ and it satisfies $f(mn) = f(m)f(n)$  whenever $(m, n) = 1$. We define {\em M\"{o}bius function} to be the function given by 
\[
\mu(n) =\left\{
                \begin{array}{ll}
                  \quad1 \quad\quad \mbox{if} \,\,\, n=1,\\
                  (-1)^k \quad \mbox{if $n=p_{i_1} p_{i_2} \cdots p_{i_k}$ are distinct primes},\\
                  \quad 0 \quad\quad \mbox{otherwise},\\
                \end{array}
              \right.
\]
or equivalently, the multiplicative function defined by $\mu(p)=-1$ and $\mu(p^k)=0$ if $k>1$ for all primes $p$. The partial sum of $\mu(n)$ function not exceeding $x$ can be defined by 
\[
M(x):=  \sum_{n \leq x} \mu(n)
\]

Several asymptotic formulas have been studied to be equivalent to the PNT by some scholars. For example, H. von Mangoldt 1897  \cite{M} proved that knowing the PNT, it is easy to obtain $\sum_{n=1}^{\infty} \frac{\mu(n)}{n} =0$ with same elementary steps. However, E. Landau 1909 showed in \cite{La} the converse of von Mangoldt's result also holds. Another equivalent of the PNT, attributed to E. Landau \cite{La1}, by $M(x)=o(x)$.

In 1912, \mbox{J. E. Littlewood} showed in \cite{Li} that the Riemann hypothesis (RH) is \mbox{equivalent} to the following evaluation
\beq \label{equation 777aa}
\begin{split} 
M(x)=  \sum_{n \leq x} \mu(n) = O(x^{\frac{1}{2}+\ep}) \quad \mbox{for all} \,\, \ep >0.
\end{split}
\eeq 
This result have been improved by some scholars (see \cite{La}, \cite{T1} and \cite{MM}).
K. \mbox{Soundararajan} 2009 \cite{S} later improved it to be
\beq \label{equation 777aax}
\begin{split} 
O\big(x^{\frac{1}{2} }\exp\big({{( \log x)^{\frac{1}{2}}} ( \log\log x)^{14}} \big) \big).
\end{split}
\eeq 
M. Balazard and A. de Roton \cite{BR} have slightly improved this bound by using  \mbox{a similar} approach as  K. \mbox{Soundararajan}. They replaced 14 by $\frac{5}{2}+\ep$ in (\ref{equation 777aax}).
The best possible bound was conjectured by S. M. Gonek (see N. Ng \cite{N}) to be 
\[
 M(x)=O\Big(x^{\frac{1}{2} } ( \log \log \log x)^{\frac{5}{4}} \Big).
\]
That is, conjecturally, one cannot get $M(x)$ to be $o\big(x^{\frac{1}{2} } ( \log \log \log x)^{\frac{5}{4}} \big)$.

It is also well-known that $M(x)$ are $\Omega(\sqrt{x})$ since there are zeros of the Riemann zeta function $\zeta$ on the line $\Re s=\frac{1}{2}$ (see for example \cite{TH}). 

\section{$CMO$ functions}
In this section, we introduce a class of functions which has been defined and \mbox{studied} by \mbox{J.-P. Kahane} and  E. Sa\"{i}as \cite{KS1}, called $CMO$ functions. These are completely \mbox{multiplicative} $f$ for which $\sum_{n=1}^\infty f(n) = 0$; $i.e.$
\[
\mbox{$f(mn) = f(m)f(n)$ \, whenever $(m, n) = 1$ \,\, and} \,\,\,\,\,\, \sum_{n=1}^\infty f(n) = 0.
\]

One of their aims was to find and give necessary and/or sufficient conditions on $f(p)$ for $f$ being a $CMO$ function. They also provided some examples of these functions. For instance, they discussed various examples of $CMO$ functions including $f(n)=\frac {\lambda(n)}{n}$, where $\lambda(n)$ is the Liouville function and $f(n)=\frac {\chi(n)}{n^\alpha}$, where $\chi$  is a non-principal Dirichlet character and $\alpha$ is a zero of $L_\chi$ with $\Re \alpha>0$.

This study drove them to think the question of how quickly partial sums of $CMO$ functions can tend to zero. They proposed that it is always $\Omega(\frac{1}{\sqrt{x}})$ and the Generalised Riemann Hypothesis -  Riemann Hypothesis (GRH-RH) would follow if their suggestion is true. This is because if GRH-RH is false then there is $\alpha$ which is a zero of $L_\chi$ with $\Re \alpha > \frac{1}{2}$ which means $\sum_{n \leq x} \frac{\chi(n)}{n^\alpha}$ is not $\Omega(\frac{1}{\sqrt{x}})$. This suggestion is incredibly difficult to prove, but it might be easier to disprove; $i.e.$, to find examples such that 
\begin{align} \label{Kahane}
  \sum_{n \leq x} f(n) = O\Big(\frac{1}{x^c}\Big) \,\,\mbox{for some} \,\, c>\frac{1}{2}. 
\end{align}

They did not find any, so in order to find example for which (\ref{Kahane}) is true we attempt to look for examples in the generalisation of $CMO$ functions.

\section{$MO$ functions}
\quad \,\, In this section,  we introduce new functions which are a natural generalisation of $CMO$ functions. We extend the notion of  $CMO$ to multiplicative functions and shall call them $MO$ functions. We would like to see how much the theory of $CMO$ \mbox{functions} can be generalised here. To help motivate our enquiries we consider \mbox{examples} of such \mbox{functions} and properties thereof. For example,  let $f$ be a $MO$ function and $g$ \mbox{a multiplicative} function ``close" to $f$.  We shall show that $g$ is also an $MO$ \mbox{function} \mbox{under} some extra condition on $f$. We can also ask a similar question of  Kahane and  Sa\"{i}as how quickly the partial sum of $MO$ functions up to and including $x$; \mbox{($i.e.$  $\sum_{n\leq x}f(n))$} can tend to zero. We define these functions as follows:


\begin{defn} \label{MO function}
An arithmetical function $f: \N \longrightarrow \C$ is called an $MO$  function if it is multiplicative and satisfies 
\[
(i) \,\, \sum_{n=1}^\infty f(n) = 0 \quad  \mbox{and} \quad  (ii) \,\, \sum_{k=0}^\infty f(p^k)  \ne 0 \,\,  \mbox{for all} \,\,  p \in \P. 
\]
\end{defn}
 The extra condition ($ii$) says the series converges but not to zero. This is needed to avoid trivial examples. For instance, let $f(1)=1$, $f(2)=-1$ and $f(n)=0$ for all $n>2$. Then $\sum_{n=1}^{\infty} f(n)= 0$  but $\sum_{k=1}^{\infty}f(2^k) = f(1)+f(2)+f(4)+ \cdots = 0$, and so does not satisfy the extra condition.

\subsection{Examples}
\quad \,\, Like $CMO$ functions which have been studied by Kahane and  Sa\"{i}as \cite{KS1},  $MO$ \mbox{functions} are not so easy to find since these need to be conditionally convergent (as we shall see in Proposition \ref{Proposition}). To help the readers understanding we give three \mbox{examples} of MO functions. The first is based on the M\"{o}bius function, the
second on the Dirichlet eta function, which corresponds to the case $k = 2$ in the third example.

\begin{exmp} \label{Mobius function}
The function $\frac{\mu(n)}{n}$ is an $MO$ function since:     

\begin{enumerate}  [label=(\roman*)] 
\item it is clear that $\frac{\mu(n)}{n}$ is a multiplicative function; 

\item it is well-known that  $\sum_{n=1}^\infty \frac {\mu(n)}{n} = 0$ (see for example \cite{A}); 
\item $\sum_{k=0}^\infty \frac{\mu(p^k)}{p^k} = 1-\frac{1}{p} \ne 0$ for all $ p \in \P$. 
\end{enumerate}

\end{exmp}

\begin{exmp}  \label{f_alpha(n) example of MO function}
Consider $\frac{(-1)^{n-1}}{n^{\alpha}}$ which is multiplicative. For which values of $\alpha \in \C$ with $\Re\alpha >0$ is this an $MO$ function?

\begin{enumerate}  [label=(\roman*)] 

\item  The series $\sum_{n=1}^\infty \frac{(-1)^{n-1}}{n^{\alpha}}$ converges for $\Re\alpha >0$  since \mbox{$A(x):=\sum_{n\leq x} (-1)^{n-1} = O(1)$.} Therefore, $0 \leq A(x) \leq 1$ and using Abel summation, we have

\beq 
\begin{split}      
 \sum_{n\leq x} \frac{(-1)^{n-1}}{n^{\alpha}} &= \frac{A(x)}{x^{\alpha}} + \alpha \int_{1}^{x}  \frac{A(t) }{t^{\alpha+1}} dt    \nonumber
\\
&= O\Big(\frac{1}{x^{\Re\alpha}}\Big) +   \alpha \int_{1}^{\infty}\frac{A(t) }{t^{\alpha+1}} dt   -  \alpha \int_{x}^{\infty}\frac{O(1) }{t^{\alpha+1}} dt  = C_\alpha+O\Big(\frac{1}{x^{\Re\alpha}} \Big), \quad 
\end{split}
\eeq

In particular, $\sum_{n=1}^\infty \frac{(-1)^{n-1}}{n^{\alpha}}$ converges. Now, for $\Re \alpha >0$, we have
\beq \label{zero}
\begin{split}
\sum_{n=1}^\infty \frac{(-1)^{n-1}}{n^{\alpha}} = (1-2^{1-\alpha})\zeta(\alpha). 
\end{split}
\eeq
This is zero if and only if $2^{\alpha} = 2$ or $\zeta(\alpha) = 0$ (for $\alpha =1$, the sum on the left \mbox{of (\ref{zero})} is not zero).
\item It remains to establish for which values of $\alpha$ that $\sum_{k=0}^\infty \frac{(-1)^{p^k-1}}{p^{\alpha k}}  \ne 0$ for all $ p \in \P$. 

If $p=2$, then 
\[
\sum_{k=0}^\infty \frac{(-1)^{2^k-1}}{2^{ k \alpha}} = 1-  \sum_{k=1}^\infty  \frac{1}{2^{\alpha k}} = \frac{2^\alpha - 2}{2^\alpha-1}.
\]
This is non-zero if  and only if  $2^{\alpha} \neq 2$; ($i.e.$ For $\frac{(-1)^{n-1}}{n^{\alpha}}$ to be $MO$ we therefore need $2^{\alpha} \neq 2$). Now if $p\geq 3$, then
\[
\sum_{k=0}^\infty \frac{(-1)^{p^k-1}}{p^{ k \alpha}} = \sum_{k=0}^\infty  \frac{1}{p^{\alpha k}} = \frac{1}{1-\frac{1}{p^\alpha}}. 
\]
This is non-zero for any $\alpha \,\,\mbox{with} \,\, \Re\alpha>0$.
\end{enumerate}

We see that $\frac{(-1)^{n-1}}{n^{\alpha}}$ is not an $MO$ function if $2^{\alpha} = 2$ since (ii) does not hold. Therefore we conclude that $\frac{(-1)^{n-1}}{n^{\alpha}}$ is an $MO$ function if and only if $\Re\alpha>0$ and $\zeta(\alpha) = 0$ since (i) and (ii) hold. 

 Furthermore, if $\zeta(\alpha) = 0$ with $\Re \alpha > 0$, then
\[
 \sum_{n \leq x}  \frac{(-1)^{n-1}}{n^{\alpha}}  =  O\Big(\frac{1}{x^{\Re \alpha}} \Big).
\]
\end{exmp}

This example can be generalised as follows:

\begin{exmp} \label{g_k(n) example of MO function}
 Define $g_k(n)$ as follows:   
\begin{align}
g_k(n) :=\left\{
                \begin{array}{ll}
                  1-k  &\mbox{if} \,\, k \,\, \mbox{divides} \,\,n,\\
                  1  &\mbox{if} \,\,k \,\, \mbox{does not divide}\,\, n.\\
                \end{array} \nonumber
              \right.
\end{align}  

We ask for which positive integer $k > 1$ and  $\alpha$ with  $\Re\alpha >0$ is the function $\frac{g_k(n)}{n^{\alpha}}$  $MO$? 
When $k=2$ we get Example \ref{f_alpha(n) example of MO function}.    

\begin{enumerate}  [label=(\roman*)] 
\item We wish to find all $k$ for which $g_k(n)$ is  a multiplicative function as follows: 
If $m= n = 1$, then $g_k (m)g_k(n)=  g_k(mn)$.
 Now if $k$ divides $mn$, then we have four cases as follows: Assume $(m,n)=1$.

\begin{enumerate}

\item If $k$ divides both $n$ and $m$, then $(m,n) \!\neq \!1$. Hence we cannot have $k$ dividing both $m,n$ since we need $(m,n) = 1$.
\item If $k$ does not  divide $n$ and $k$ divides $m$, then $g_k (m)g_k(n)= (1-k)(1) = 1-k = g_k(mn)$. 

or vice versa
\item If $k$  does not  divide $m$ and $k$ divides $n$, then $g_k (m)g_k(n)= (1)(1-k)= 1-k = g_k(mn)$. 

\item If $k$ does not divide both $n$ and $m$, then we have two cases: 

\begin{enumerate}

\item  If $k$ is not a prime power; ($i.e.$ $k=p_1^{a_1}\!\cdot p_2^{a_2} \cdots p_i^{a_i}, \,\, \mbox{where}\,\, i \geq 2 \,\mbox{and}$ \mbox{$a_i \geq 1$}). Then, with $m=p_1^{a_1}$ and $n=p_2^{a_2} \cdots p_i^{a_i}$ such that $(m,n)=1$, we have $g_k (m)g_k(n)= (1)(1) \neq (1-k) =  g_k(mn)$.   

\item If $k$ is a prime power; ($i.e.$ $k=p^r$). Then at least one of $m$ or $n$ is not \mbox{a multiple} of $p$ while the other is  ($i.e.$ $p$ does not divide $m$, then $p^r$ divides $n$ or  $p$ does not divide $n$, then $p^r$ divides $m$) and 
$g_k (m)g_k(n)= (1)(1-k) = (1-k) =  g_k(mn)$ or $g_k (m)g_k(n)= (1-k)(1) = (1-k) =  g_k(mn)$. 

\end{enumerate}

\end{enumerate}\
However,  if \,$k$\, does not divide $mn$, then \,$k$ \,does not  divide both\, $m$ \,and \,$n$,\, and \,$g_k (m)g_k(n)= (1)(1) = 1 = g_k(mn)$.

Thus $g_k(n)$ is multiplicative function if and only if $k$ is a prime power. 

\item  The series $\sum_{n=1}^\infty \frac{g_k(n)}{n^{\alpha}}$ converges for $\Re\alpha >0$  since
\beq 
\begin{split}
A(x):=\sum_{n\leq x} g_k(n)&= \sum_{m=1}^{N} \,\,\sum_{n=(m-1)k+1}^{mk} g_k(n) + \sum_{n=Nk+1}^{x} g_k(n) = 0+ \sum_{n=Nk+1}^{x} g_k(n) \nonumber
\\
&= g_k(Nk+1) +g_k(Nk+2)+ \cdots + g_k(x),  \quad \mbox{where} \,\, N=\left\lfloor{\frac{x}{k}}  \right \rfloor
\\
&\leq  k-1=O(1).
\end{split}
\eeq
Thus $0 \leq A(x) \leq k-1$ and using Abel summation, we have
\beq 
\begin{split}      
 \sum_{n\leq x} \frac{g_k(n)}{n^{\alpha}} &= \frac{A(x)}{x^{\alpha}} + \alpha \int_{1}^{x}  \frac{A(t) }{t^{\alpha+1}} dt    \nonumber
\\
&= O\Big(\frac{1}{x^{\Re\alpha}}\Big) +   \alpha \int_{1}^{\infty}\frac{A(t) }{t^{\alpha+1}} dt   -  \alpha \int_{x}^{\infty}\frac{O(1) }{t^{\alpha+1}} dt      
\\
&= C_\alpha+O\Big(\frac{1}{x^{\Re\alpha}} \Big), \quad 
\end{split}
\eeq
where $C_\alpha$ is a constant, as in Example \ref{f_alpha(n) example of MO function}.
In particular, for $\Re \alpha>0$, $\sum_{n=1}^\infty \frac{g_k(n)}{n^{\alpha}}$ converges.

Now, for $\Re \alpha>1$, we have
\beq 
\begin{split}
\sum_{n=1}^\infty \frac{g_k(n)}{n^{\alpha}} &=\sum_{n=1}^\infty \frac{1}{n^\alpha}- \sum_{n=1}^\infty \frac{k}{(kn)^\alpha} =(1-k^{1-\alpha})\zeta(\alpha).  \nonumber
\\
\mbox{Thus} \,\, \sum_{n=1}^\infty \frac{g_k(n)}{n^{\alpha}} &=  C_\alpha  =(1-k^{1-\alpha})\zeta(\alpha) \quad \mbox{for} \,\, \Re \alpha>0 \,\, \mbox{by analytic continuation.} 
\\
\mbox{Also}, \,\,  \sum_{n=1}^\infty \frac{g_k(n)}{n^{\alpha}}  &=0 \,\, \mbox{if and only if} \,\,k^{\alpha} = k\,\, \mbox{or} \,\,\zeta(\alpha) = 0.
\end{split}
\eeq
\item It remains to get all $k$ and $\alpha$ for which $\sum_{m=0}^\infty \frac{g_k(p^m)}{p^{m \alpha}}  \ne 0 $ for all $ p \in \P$. Let  $k=p_0^r$, $p_0$ a prime number.

 If $p_0\ne p$, then $g_k(p^m)=1$ for all $m \geq 0$. Hence
\[
\sum_{m=0}^\infty \frac{g_k(p^m)}{p^{ m\alpha }} = \sum_{m=0}^\infty  \frac{1}{p^{ \alpha m}} = \frac{1}{1-\frac{1}{p^\alpha}}.
\]
This is non-zero for any $\alpha \,\,\mbox{with} \,\Re \alpha>0$. Now if $p_0=p$, then
\beq 
\begin{split}
\sum_{m=0}^\infty \frac{g_k(p^m)}{p^{m \alpha}} &= \sum_{m=0}^{r-1}  \frac{g_k(p^m)}{p^{\alpha m}} +\sum_{m=r}^\infty  \frac{g_k(p^m)}{p^{\alpha m}} =    \sum_{m=0}^{r-1}  \frac{1}{p^{\alpha m}} +(1-k)\sum_{m=r}^\infty  \frac{1}{p^{\alpha m}}
\\
&= \frac{1-p^{-\alpha m}}{1-p^{-\alpha}} + (1-k)\frac{1}{p^{r\alpha}}\frac{1}{1-p^{-\alpha}}=  \frac{p^{r\alpha}-k}{p^{(r-1)\alpha}(p^{\alpha} -1)} = \frac{k^{\alpha}-k}{k^{\alpha}(1-p^{-\alpha})}. \nonumber
\end{split}
\eeq
This is non-zero if  and only if $k^{\alpha} \neq k$; ($i.e.$, For $\frac{g_k(n)}{n^{\alpha}}$ to be $MO$ we therefore need $k^{\alpha} \neq k$). 
\end{enumerate}

We see that $\frac{g_k(n)}{n^{\alpha}}$ is not an $MO$ function if $k^{\alpha} = k$ since (iii) fails. Therefore, we conclude that $\frac{g_k(n)}{n^{\alpha}}$ is an $MO$ function if and only if $k$ is a prime power, $\Re \alpha > 0$ and $\zeta(\alpha) = 0$  since (i), (ii) and (iii)  hold.

 Furthermore, if $\zeta(\alpha) = 0$ with $\Re \alpha > 0$, then
\[
 \sum_{n \leq x}  \frac{g_k(n)}{n^{\alpha}}  =  O\Big(\frac{1}{x^{\Re \alpha}} \Big).
\]
\end{exmp}

\subsection{Some properties of $MO$ functions}
\quad \,\, In this section, we establish some preliminary properties of $MO$ functions. We shall first need the following result in the course of our discussion.

\begin{pro} \label{Proposition a110}
Let $f$ be a multiplicative function. Then $\sum_{n =1}^{\infty} |f(n)| $ converges, so that $f$ is absolutely convergent, if and only if $\sum_{p} \sum_{k =1}^{\infty}  |f(p^k)|$  converges.
\end{pro}

\begin{prf}
Trivially, the series $\sum_{p} \sum_{k =1}^{\infty}  |f(p^k)|$ converges if $\sum_{n =1}^{\infty} |f(n)|$ converges.

Now suppose  $\sum_{p} \sum_{k =1}^{\infty}  |f(p^k)|$  converges. 
It is follows that
\beq 
\begin{split}   
\prod_{p}\Big(1+ \sum_{k =1}^{\infty} |f(p^k)| \Big) = \prod_{p}\Big( \sum_{k =0}^{\infty} |f(p^k)| \Big)\,\, \mbox{ converges}. \nonumber
\end{split}
\eeq
But the right hand side is at least $\prod_{p \leq x} \Big\{ \sum_{k =0}^{\infty} |f(p^k)| \Big\}$. Therefore, by the proof of  \mbox{Theorem 11.6} of \cite{A} for any $x$, we have
\beq 
\begin{split} 
\prod_{p \leq x} \Big\{ \sum_{k =0}^{\infty} |f(p^k)| \Big\} =\!\!\! \sum_{\tiny\begin{array}{c}  n \in \N \\  p|n \,\,\& \,\,p \leq x \end{array}}\!\!\! |f(n)| \geq \sum_{n \leq x} |f(n)|. \nonumber
\end{split}
\eeq
Hence $\sum_{n =1}^{\infty} |f(n)| \,\,  \mbox{converges, so that $f$ is absolutely convergent.} $
\bo
\end{prf}

\begin{pro} 
If $f$ is a $CMO$ function, then  $f$ is an $MO$ function; $(i.e \,\,CMO \subset MO) $.
\end{pro}

\begin{prf}

It is clear that $f$ is multiplicative and $\sum_{n=1}^{\infty} f(n) = 0$. It remains to show that $\sum_{k=0}^{\infty}f(p^k) \neq 0$ for all $ p \in \P$.
Now since $f$ is completely multiplicative, then $f(p^k) = f(p)^k$. Therefore
 \beq 
\begin{split}  
\sum_{k=0}^{\infty}f(p^k) &= \sum_{k=0}^{\infty}f(p)^k  = \frac{1}{1-f(p)} \neq 0.  \nonumber
\end{split}
\eeq
 \\
This series converges since $|f(p)| < 1$. Hence, by Definition \ref{MO function}, $f $ is an $MO$ function.
\bo
\end{prf}

\begin{pro}  \label{Proposition}
Let \,$f$ be an \,$MO$ \,function. Then \,$\sum_{n=1}^{\infty} |f(n)|$ \,diverges. Indeed  $\sum_{p} \sum_{k=1}^\infty |f(p^k)|$ diverges.
\end{pro}

\begin{prf}

Let us assume that the statement is false, so that
\beq \label{equation 31}
\begin{split}
\sum_{n=1}^{\infty} |f(n)| \,\, \mbox{ converges.} \nonumber
\end{split}
\eeq
Then, by multiplicative property,    
\[
\sum_{n=1}^{\infty} f(n) = \prod_{p} \sum_{k=0}^{\infty} f(p^k) \neq 0 \,\, \mbox{since} \,\,  \sum_{k=0}^{\infty} f(p^k) \neq 0.
\]
Yielding a contradiction since $f$ is an $MO$ function and  hence 
\beq \label{equation oo12234}
\begin{split}
\sum_{n=1}^{\infty} |f(n)|  \,\,\, \mbox{diverges.}  \nonumber
\end{split}
\eeq
 
Furthermore, Proposition \ref{Proposition a110} gives $\sum_{p} \sum_{k=1}^\infty |f(p^k)|$ diverges, as required. 
\bo
\end{prf}

\subsubsection{Partial sums of $MO$ functions}
\quad \,\, We know that the partial sum of an $MO$ function not exceeding $x$ tends to zero when $x$ tends to infinity. A question raised by Kahane and  Sa\"{i}as \cite{KS1} regarding $CMO$ functions is: can one show, given $g(x)$, that there exist a $CMO$ function $f$ with

\[
\sum_{n \leq x} f(n) = \Omega(g(x))?
\]
We are not considering this question, but we are interested in a related question which is: how small can we make $g(x)$, so that the above is true for all $MO$ functions $f$? This question motivates the following propositions:

\begin{pro} \label{Proposition lower bounds of MO function}
If $f$ is an $MO$ function, then $$\sum_{n \leq x} f(n) = \Omega\Big( \frac{1}{x \log x} \Big).$$
\end{pro}

\begin{prf}

Let us assume that the statement is false, so that
\[
 \sum_{n \leq x} f(n) = O\Big( \frac{1}{x \log x} \Big).
\]

We know that for $n \in \N$,
$$ f(n) =\sum_{m \leq n} f(m)- \sum_{m < n} f(m) = O\Big( \frac{1}{n \log n} \Big).$$
Hence
$$ f(p^k) = O\Big( \frac{1}{p^k \log p^k} \Big).$$

Now it follows that $\sum_{p} \sum_{k =1}^{\infty} |f(p^k)|$ converges since
\beq 
\begin{split}   
\sum_{p} \sum_{k =1}^{\infty} \frac{1}{p^{k}\log p^k} &\leq \sum_{p}\sum_{k =1}^{\infty} \frac{1}{p^{k} \log p} \quad (\mbox{since}\,\, \log p^{k} \geq \log p ) \nonumber
\\
&= \sum_{p} \frac{1}{\log p} \, \sum_{k =1}^{\infty}  \frac{1}{p^{k}}
\\
&= \sum_{p} \frac{1}{(p-1)\log p}   \,\, \mbox{ converges (since $p_n \log p_n \sim n(\log n)^2$).} 
\end{split}
\eeq
Thus
\beq 
\begin{split}   
 \sum_{p} \sum_{k =1}^{\infty} \frac{1}{p^{k}\log p^k} \,\,\mbox{converges.} \nonumber
\end{split}
\eeq

Hence, by Proposition \ref{Proposition a110}, $\sum_{n=1}^{\infty} |f(n)|$ converges. However, by Proposition \ref{Proposition}, we have a contradiction, and so it follows that
\[
\sum_{n \leq x} f(n) = \Omega\Big( \frac{1}{x \log x} \Big).
\]
\bo
\end{prf}

\begin{rmk} 
Similarly, if $f$ is an $MO$ function, then
\[
\sum_{n \leq x} f(n)= \Omega\big( \frac{1}{x (\log x)^{\ep}} \big) \quad \mbox{for all} \,\, \ep > 0.
\]

We can improve Proposition \ref{Proposition lower bounds of MO function} using the fact that $\sum_{p} \frac{1}{p(\log \log p)^2}$ converges.
\end{rmk}

\begin{pro} \label{ Proposition 49}
If $f$ is an $MO$ function, then $$\sum_{n \leq x} f(n) = \Omega\Big( \frac{1}{x (\log \log x)^2} \Big).$$
\end{pro}

\begin{prf}

Let us assume that the statement is false, so that
\[
 \sum_{n \leq x} f(n) = O\Big( \frac{1}{x (\log \log x)^2} \Big).
\]

We know that for $n \in \N$,
$$ f(n) =\sum_{m \leq n} f(m)- \sum_{m < n} f(m) = O\Big( \frac{1}{n (\log\log n)^2} \Big).$$
Hence
$$ f(p^k) = O\Big( \frac{1}{p^k (\log \log p^k)^2} \Big).$$

Now it follows that $\sum_{p} \sum_{k =1}^{\infty} |f(p^k)|$ converges
since
\beq 
\begin{split}   
\sum_{p \geq 3} \sum_{k =1}^{\infty} \frac{1}{p^{k}(\log \log p^k)^2} &\leq \!\sum_{p \geq 3} \sum_{k =1}^{\infty} \frac{1}{p^{k}(\log \log p)^2} \quad (\mbox{since}\,\, (\log \log p^{k})^2 \geq (\log \log p)^2 ) \nonumber
\\
&= \!\sum_{p \geq 3} \frac{1}{(\log \log p)^2} \, \sum_{k =1}^{\infty}  \frac{1}{p^{k}}
\\
&=\! \sum_{p \geq 3} \frac{1}{(p-1)(\log \log p)^2}   \,\mbox{converges (since $(\log\log p_n)^2 \!\sim\! (\log \log n)^2$).} 
\end{split}
\eeq
For $p=2$, 
\beq 
\begin{split}   
\sum_{k =1}^{\infty} \frac{1}{2^{k}(\log \log 2^k)^2} \leq \frac{1}{2(\log \log 2)^2 }+ \frac{1}{(\log \log4)^2}\sum_{k \geq 2 } \frac{1}{2^{k}} \,\,\,\,\mbox{converges.} \nonumber
\end{split}
\eeq
Thus
\beq 
\begin{split}   
 \sum_{p} \sum_{k =1}^{\infty} \frac{1}{p^{k}(\log \log p^k)^2} \,\,\mbox{converges.} \nonumber
\end{split}
\eeq

Hence,  by  Proposition \ref{Proposition a110}, $\sum_{n=1}^{\infty} |f(n)|$ converges. However, by Proposition \ref{Proposition}, we have a contradiction, and so it follows that 
\[
\sum_{n \leq x} f(n) = \Omega\Big( \frac{1}{x (\log \log x)^2} \Big).
\]
\bo
\end{prf}

\begin{rmk} \label{Remark 50}
Similarly, if $f$ is an $MO$ function, then
\beq \label{Similarly}
\begin{split}  
\sum_{n \leq x} f(n)= \Omega\big( \frac{1}{x (\log \log x)^{1+\ep}} \big) \quad \mbox{for all} \,\, \ep > 0.
\end{split}
\eeq
\end{rmk}

 Kahane and  Sa\"{i}as \cite{KS1} have shown that if $f$ is a $CMO$ function, then 
\[
\sum_{n \leq x} f(n)= \Omega\Big( \frac{1}{x} \Big)
\]
by using a deep result of D. Koukoulopoulos in \cite{K}. We attempted to improve (\ref{Similarly}) to $\Omega(\frac{1}{x})$ as with the work of Kahane and Sa\"{i}as, but the question is still open.

\subsubsection{Closeness relation between two multiplicative functions}
\quad \,\, Let $\mathcal{M} := \{ f:\N \longrightarrow \C \,\, \mbox{multiplicative}\}$, and let us define an {\em (extended) metric} on $\mathcal{M}$ to be the distance function $$D(f,g):=\sum_{p}\sum_{k=0}^{\infty} |g(p^k)- f(p^k)|.$$
Then $ \mathcal{M}$ is an {\em extended metric space} since $D(f,g)$ can attain the value $\infty$. It is \mbox{straightforward} to check for all $f,g,h \in \mathcal{M}$ 
\begin{enumerate}[label=(\roman*)]
\item $D(f,g)=0$ if and only if $f=g$,
\item $D(f,g) = D(g,f)$,
\item $D(f,h) \leq D(f,g) +D(g,h)$,
\end{enumerate}
hold. We aim to extend Theorem 3 of Kahane and Sa\"{i}as in \cite{KS1} by showing that if $f$ is an $MO$ function and $g$ is a multiplicative function ``close" to $f$, ($i.e.$ $g$ has finite distance from $f$), then $g$ is also an $MO$ function. We can do this under an extra condition on $f$, as the following theorem shows.

\begin{thm} \label{ttww}
Let $f$ be an $MO$ function for which 
\beq \label{e}
\begin{split}   
 \Big|\sum_{k=0}^{\infty} \frac{f(p^k)}{p^{ks}}\Big| \geq a \quad \mbox{for some} \,\, a>0, \,\, \mbox{for all} \,\,p \,\, \mbox{and all} \,\, \Re s\geq 0, 
\end{split}
\eeq 
and let $g$ be a multiplicative function such that $D(f,g)$ is finite and
\beq  \label{eee}
\begin{split}   
 \sum_{k=0}^{\infty} g(p^k) \ne 0 \quad \mbox{for all} \,\,p.
\end{split}
\eeq
Then $g$ is an $MO$ function.
\end{thm}

\begin{prf}

Let $F(s):=\sum_{n=1}^{\infty} \frac{f(n)}{n^s}$ and $G(s):=\sum_{n=1}^{\infty} \frac{g(n)}{n^s}$. Then the series for $F(s)$ is absolutely convergent for $\Re s>1$ and it is convergent for $\Re s >0$ and $s=0$ since $\sum_{n=1}^{\infty} f(n) = 0$. We note that the assumption $D(f,g)$ is finite and the fact that  $f$ is an  $MO$ function imply $|g(p^k)| \longrightarrow 0$ as $p^k \longrightarrow \infty$. Then, by Theorem 316 of \cite{HW}, $g(n) \longrightarrow 0$ as $n \longrightarrow \infty$. Therefore the series for $G(s)$ converges for $\Re s >1$ since $g$ is bounded.  Therefore $F(s)$ and $G(s)$ can be written as follows:
\[
 F(s)=\prod_{p}\sum_{k=0}^{\infty} \frac{f(p^k)}{p^{ks}} \quad \mbox{and} \quad G(s)=\prod_{p}\sum_{k=0}^{\infty} \frac{g(p^k)}{p^{ks}}  \quad \Re s >1.
\]
Now
\beq 
\begin{split}  
H(s):=\prod_{p} \Bigg(\frac{\sum_{k=0}^{\infty} \frac{g(p^k)}{p^{ks}}}{\sum_{k=0}^{\infty} \frac{f(p^k)}{p^{ks}}}\Bigg) = \prod_{p} \Bigg(1+\frac{\sum_{k=0}^{\infty} \frac{g(p^k)-f(p^k)}{p^{ks}}}{\sum_{k=0}^{\infty} \frac{f(p^k)}{p^{ks}}} \Bigg)    \nonumber
\end{split}
\eeq
converges absolutely for $\Re s \geq 0$ if and only if 
\beq \label{22222}
\begin{split}
\sum_{p}  \frac{\big|\sum_{k=0}^{\infty} \frac{ g(p^k)-f(p^k)}{p^{ks}}\big|}{\big|\sum_{k=0}^{\infty} \frac{f(p^k)}{p^{ks}}\big|} 
\end{split}
\eeq
converges for $\Re s \geq 0$. But
\beq 
\begin{split}
\sum_{p}  \frac{\big|\sum_{k=0}^{\infty} \frac{ g(p^k)-f(p^k)}{p^{ks}}\big|}{\big|\sum_{k=0}^{\infty} \frac{f(p^k)}{p^{ks}}\big|} \leq \frac{1}{a} \sum_{p}\sum_{k=0}^{\infty}  |g(p^k)- f(p^k)| \nonumber
\end{split}
\eeq
by (\ref{e}) so, since $D(f,g)$ is finite, (\ref{22222}) converges for $\Re s \geq 0$ and $H(s)$ converges absolutely to holomorphic function for $\Re s > 0$. However, $H(s)=(G/F)(s)$ for $\Re s >1$ then $G(s)=F(s)H(s)$, where  the series for $F(s)$ converges for $\Re s> 0$ and $s=0$ since $f$ is an $MO$ function, and $H(s)$ converges absolutely for $\Re s\geq 0$. Therefore $G(s)$ converges for $\Re s> 0$ and $s=0$  using the extension of Theorem 1.2 of Chapter II.1. \cite{T}. Thus we have $G(0)=F(0)H(0)=0$. Hence, by assumption (\ref{eee}) and $G(0)=0$, $g$  is an $MO$ function.
\bo
\end{prf}

The proof of Theorem \ref{ttww} also gives the following result.

\

\begin{cor} 
Let $f$ and $g$ both be  multiplicative functions such that $D(f,g)$ is finite and satisfies 
\beq \label{eeee}
\begin{split}   
 \Big|\sum_{k=0}^{\infty} \frac{f(p^k)}{p^{ks}}\Big| \geq a \quad \mbox{for some} \,\, a>0 \,\, \mbox{and all} \,\, \Re s\geq 0,  \nonumber
\end{split}
\eeq 
\beq \label{eeeee}
\begin{split}   
 \Big|\sum_{k=0}^{\infty} \frac{g(p^k)}{p^{ks}}\Big| \geq b \quad \mbox{for some} \,\, b>0 \,\, \mbox{and all} \,\, \Re s\geq 0.  \nonumber
\end{split}
\eeq 

Then the following two assertions are equivalent:
\beq \label{eeeeeee}
\begin{split}   
\sum_{n=1}^{\infty} f(n)=0 \quad \mbox{and} \quad \sum_{n=1}^{\infty} g(n)=0. \nonumber
\end{split}
\eeq 
\end{cor} 

\subsection{Open problems}
\begin{enumerate}[label=(\roman*)] 
\item  Let $f$ be an $MO$ function. Can we show that 
\[
  \sum_{n \leq x } f(n) = \Omega\Big(\frac{1}{x}\Big) \,\, ?
\]

\item As pointed out earlier Kahane and Sa\"{i}as  suggested that for all $CMO$ \mbox{functions,} one \,has $\sum_{ n \leq x} f(n) =  \Omega\big(\frac{1}{\sqrt{x}} \big)$. As \,also\, mentioned, GRH-RH (Generalised \mbox{Riemann} Hypothesis-Riemann Hypothesis) would follow if their suggestion is correct.

In Example \ref{Mobius function}, it is known that $\sum_{n\leq x} \mu(n)=\Omega(\sqrt{x})$ since there are zeros of the Riemann zeta function $\zeta$ on the line $\Re s=\frac{1}{2}$ (see \cite{T}). Thus,  by Abel summation, 
\[
\sum_{n\leq x} \frac{\mu(n)}{n} =   \Omega\Big(\frac{1}{\sqrt{x}} \Big). 
\] 

However, for $\sum_{n \leq x}  \frac{(-1)^{n-1}}{n^{\alpha}}$ and $\sum_{n \leq x}  \frac{g_k(n)}{n^{\alpha}}$ to converge to zero  in Examples \ref{f_alpha(n) example of MO function} and \ref{g_k(n) example of MO function}, it is necessary that $\alpha$ be a zero of $\sum_{n =1}^\infty  \frac{(-1)^{n-1}}{n^{s}}$ and $\sum_{n =1}^\infty  \frac{g_k(n)}{n^{s}}$ with $\Re \alpha > 0$; ($i.e.$ $\zeta(\alpha)=0$). Suppose this is the case. We then have 
\[
 \sum_{n \leq x}  \frac{(-1)^{n-1}}{n^{\alpha}} = O\Big(\frac{1}{x^{\Re \alpha}} \Big) \,\,\mbox{and}\,\, \sum_{n \leq x}  \frac{g_k(n)}{n^{\alpha}}= O\Big(\frac{1}{x^{\Re \alpha}} \Big),
\]
and 
\[
 \sum_{n \leq x}  \frac{(-1)^{n-1}}{n^{\alpha}} = \Omega\Big(\frac{1}{x^{\Re \alpha}} \Big) \,\,\mbox{and}\,\, \sum_{n \leq x}  \frac{g_k(n)}{n^{\alpha}}= \Omega\Big(\frac{1}{x^{\Re \alpha}} \Big).
\]

\end{enumerate}

In our results, we have not found any examples with $\sum_{ n \leq x} f(n) =  O\big(\frac{1}{x^c} \big)$ for $c > \frac{1}{2}$. This may suggest the following conjecture. 

\begin{con}  \label{Conjecture MO} 
 For all  multiplicative function $f$  ($MO$ functions), we have 
\[   
     \sum_{ n \leq x } f(n) =  \Omega\Big(\frac{1}{\sqrt{x}} \Big).        
 \] 
\end{con}

Furthermore, the RH would follow if Conjecture \ref{Conjecture MO} were true since if RH is false then there is  $\alpha$ which is a zero of $\zeta$ with $\Re \alpha > \frac{1}{2}$ which means $\sum_{n \leq x} \frac{(-1)^{n-1}}{n^\alpha}$ and $\sum_{n \leq x}  \frac{g_k(n)}{n^{\alpha}}$  is not $\Omega(\frac{1}{\sqrt{x}})$.

{\small


\begin{thebibliography}{99}

\bibitem{A} T. M. Apostol, {\em Introduction to Analytic Number Theory}, Springer, 1976.

\bibitem{BR} M. Balazard and A. D. Roton. Notes de lecture de l'article Partial sums of the Mobius function de Kannan Soundararajan, arXiv preprint arXiv:0810.3587, 2008.

\bibitem{HW} G. H. Hardy and E. M. Wrigh, {\em An introduction to the theory of numbers}, Oxford university press, 1979.

\bibitem{K} D. Koukoulopoulos, On multiplicative functions which are small on average, { \em Geometric and Functional Analysis}, {\bf 23}(5):1569-1630, 2013.

\bibitem{KS1} J.-P. Kahane and E. Sa\"{i}as,  Fonctions compl\'{e}tement multiplicatives de somme nulle,  Expositiones Mathematicae, {\bf 35}(4):364-389, 2017.

\bibitem{Li} J. E. Littlewood, Quelques consequences de l'hypothese que la fonction $\zeta(s)$ de Riemann n'a pas de zeros dans le demi-plan Rs > 1
2, CRAS Paris, 154:263-266, 1912.

\bibitem{La} E. Landau, Uber die mobiussche funktion. Rendiconti del Circolo Matematico di Palermo (1884-1940), 48(2):277-280, 1924.

\bibitem{La1} E. Landau, Neuer beweis der gleichung $\sum_{k=0}^\infty \frac{mu(k)}{k}$, PhD thesis, Druck Der Dieterich'schen Univ.-Buchdruckerei (W. Fr. Kastner), 1899.

\bibitem{M} H. von Mangoldt, Beweis der gleichung $\sum_{k=0}^\infty \frac{mu(k)}{k}$, Proc. Royal Pruss. Acad. of Sci. of Berlin, pages 835-852, 1897.

\bibitem{MM} H. Maier and H. L. Montgomery, The sum of the Mobius function. Bulletin of the London Mathematical Society, 41(2):213-226, 2009.

\bibitem{N} N. Ng,  The distribution of the summatory function of the Mobius function. Proceedings of the London Mathematical Society, 89(2):361-389, 2004.

\bibitem{S} K. Soundararajan, Partial sums of the Mobius function. Journal fur die reine und angewandte Mathematik (Crelles Journal), 361:141-152, 2009.

\bibitem{T} G. Tenenbaum, {\em Introduction to analytic and probabilistic number theory}, volume 163. American Mathematical Society, 2015.

\bibitem{TH} E. C. Titchmarsh and D. R Heath-Brown, {\em The theory of the Riemann zeta-function}, Oxford University Press, 1986.

\bibitem{T1} E. C. Titchmarsh, A consequence of the Riemann Hypothesis, Journal of the London Mathematical Society, 1(4):24-254, 1927.

\end{thebibliography}
\end{document}